\documentclass[leqno]{article}
\usepackage{amssymb,color}

\newtheorem{theorem}{Theorem}

\newtheorem{remark}[theorem]{Remark}

\begin{document}
\title{Four-dimensional Osserman metrics with nondiagonalizable Jacobi operators\footnote{\rm Supported by project  BFM 2003-02949 (Spain)}}
\author{J. Carlos D\'{\i}az-Ramos\and Eduardo Garc\'{\i}a-R\'{\i}o\and Ram\'{o}n V\'{a}zquez-Lorenzo}
\date{}
\maketitle
\begin{abstract}
A complete description of Osserman four-manifolds whose Jacobi operators have
a nonzero double root of the minimal polynomial is given.
\end{abstract}

\vspace{0.5cm}

\begin{enumerate}
\item[] {\underline{\normalsize\bf 2000 M. S. C.:}} 53C50, 53B30
\item[] {\underline{\normalsize\bf Key words:}}
Jacobi operator, Jordan-Osserman metric,
Walker metric, self-dual and anti-self-dual Weyl tensor.
\end{enumerate}

\section{Introduction}

Geometric information about a pseudo-Riemannian manifold $(M,g)$ is essentially
encoded by the curvature tensor $R\in\otimes^4T^*M$. Hence a central problem in
differential geometry is to relate algebraic properties of the curvature tensor
to the underlying geometry of the manifold. Due to the fact that the whole curvature
tensor is so difficult to handle, the investigation usually focuses on different
objects associated to the curvature tensor, the Jacobi operator being the most
natural and widely investigated.
A pseudo-Riemannian manifold $(M,g)$ is said to be \emph{Osserman}
if the eigenvalues of the Jacobi operators are constant on the
unit pseudo-sphere bundles $S^\pm (TM)$.
Since the Ricci tensor is obtained by tracing the Jacobi operators, any
Osserman metric is Einstein and, in particular, of constant sectional curvature in
dimension $2$, $3$.

Dimension four is therefore the first non-trivial case
for consideration in trying to determine the whole class of Osserman metrics.
Moreover due to the connection between Osserman and (anti-) self-dual metrics many
special features occur like the existence of \emph{pointwise Osserman} metrics
which are not Osserman, i.e., the eigenvalues of the Jacobi operators are still
independent of the direction, but they may change from point to point.
(See for example \cite{GKV}, \cite{Gilkey1}, \cite{GSV} and the references therein
for more information on pointwise Osserman $4$-manifolds,
and \cite{BCGHV}, \cite{G-I}, \cite{G-I-2}, \cite{Ni1}, \cite{Ni2}
for the higher dimensional Riemannian and pseudo-Riemannian cases).

It has been shown in \cite{Chi} and \cite{GKVa} (see also \cite{BBG}) that any
four-dimensional Osserman metric is
locally isometric to a two-point homogeneous space if the signature is either Riemannian or Lorentzian.
However the situation is much more complicated for neutral metrics and there exist many examples  of
nonsymmetric Osserman pseudo-Riemannian manifolds of neutral signature \cite{GR-VA-VL}. Indeed besides the partial
results at \cite{BBR} for Osserman $(++--)$-metrics with diagonalizable Jacobi operators and those at
\cite{GV} for locally symmetric Osserman four-manifolds, the general problem of obtaining a complete
description of four-dimensional Osserman metrics of neutral signature remains open.

For any non-null vector $X$ in the $(++--)$-setting, the induced metric on $X^\perp$ is of Lorentzian signature, and hence
the eigenvalue-structure does not completely characterize the Jacobi operators $R_X$. The consideration
of the Jordan normal form led to introduce the so-called \emph{Jordan-Osserman} metrics as those where
the Jordan normal form of the Jacobi operators is constant on $S^\pm (TM)$. Four-dimensional Jordan-Osserman
metrics were initially investigated by Bla\v zi\'c, Bokan and Raki\'c \cite{BBR},
who considered four different
possibilities according to the behavior of the Jordan normal form of the Jacobi operators as follows:
\begin{itemize}
\item[(\emph{Ia}):] The Jacobi operators are diagonalizable, i.e.,
$$
R_X=\left(\begin{array}{ccc}
\alpha & &\\
       &\beta&\\
       & &\gamma
\end{array}\right).
$$
\item[(\emph{Ib}):] The Jacobi operators have a complex eigenvalue, i.e.,
$$
R_X=\left(\begin{array}{ccc}
\alpha &-\beta&\\
\beta&\alpha&\\
       & &\gamma
\end{array}\right).
$$
\item[(\emph{II}):] There is a double root of the minimal polynomial of $R_X$, i.e.,
$$
R_X=\left(\begin{array}{ccc}
\alpha & &\\
       &\beta&\\
       & 1&\beta
\end{array}\right).
$$
\item[(\emph{III}):] There is a triple root of the minimal polynomial of $R_X$, i.e.,
$$
R_X=\left(\begin{array}{ccc}
\alpha & &\\
    1   &\alpha&\\
       & 1&\alpha
\end{array}\right).
$$
\end{itemize}

It is shown in \cite{BBR} that Jordan-Osserman metrics with diagonalizable Jacobi operators are locally isometric
to a real, complex or paracomplex space form and also that Type Ib metrics cannot occur. Moreover, locally
symmetric Jordan-Osserman $(++--)$-metrics have diagonalizable Jacobi operators or they are isometric to some Type II
metrics with nilpotent Jacobi operators \cite{GV}.
The fact that all known examples of nonsymmetric Jordan-Osserman metrics had nilpotent Jacobi operators of degree
two or three suggested
that no other examples could exist (cf. \cite{CGM}, \cite{GKV}, \cite{Gilkey1}).
However, very recently the authors have shown the existence of a family of
Type II Jordan-Osserman metrics with non-nilpotent Jacobi operators \cite{DGV}.
The purpose of this work is to clarify the situation of Type II Jordan-Osserman metrics by proving the following

\bigskip

\noindent{\normalsize\bf Main Theorem}
\emph{Let $(M,g)$ be a four-dimensional Type II Jordan-Osserman manifold.
Then the Jacobi operators are two-step nilpotent or otherwise there
exist local coordinates $(x_1,\dots,x_4)$ where the metric is given by
\begin{equation}\label{typeii-1}
dx^1 dx^3+ dx^2 dx^4 + \sum_{i\leq j=3,4} s_{ij} dx^i dx^j
\end{equation}
for some functions $s_{ij}(x_1,\dots,x_4)$ as follows
\begin{equation}\label{typeii-2}
\begin{array}{l}
 s_{33} =\! x_1^2 \frac{\tau}{6} \!+ x_1 P \!+ x_2 Q\!
   + \frac{6}{\tau}\! \left\{ Q   ( T\! -\! U ) + V   (P\! -\! V) - 2 (Q_4\!-\!  V_3) \!\right\}\! ,\!\!\!
   \\[0.1in]
 s_{44} =\!  x_2^2 \frac{\tau}{6} \!+ x_1 S \!+ x_2 T\!
   + \frac{6}{\tau}\! \left\{ S   ( P \! - \! V ) + U   (T\! -\! U) - 2 (S_3\! -\! U_4) \!\right\}\! ,\!\!\!
   \\[0.1in]
 s_{34} =\!  x_1 x_2 \frac{\tau}{6} \!+ x_1 U \!+  x_2 V\!
   + \frac{6}{\tau}\! \left\{ - Q S + U V + T_3 - U_3 + P_4 - V_4 \!\right\}\! ,\!\!\!
\end{array}
\end{equation}
and arbitrary functions $P$, $Q$, $S$, $T$, $U$, $V$ depending only on the coordinates
$(x_3,x_4)$, where $\tau\neq 0$ denotes the scalar curvature.
}

\bigskip


The proof of the Main Theorem is based on the following facts.

\medskip
\noindent{\normalsize\bf Fact A.} \,{\rm\cite{ABBR}, \cite{GSV}}\label{th: Osserman = Einstein self-dual (anti-self-dual)}
{A four-dimensional pseudo-Riemannian manifold  is pointwise Osserman if and only if
it is Einstein self-dual (or anti-self-dual)}.
\medskip

\noindent {\normalsize\bf Fact B.} \cite[Corollary 8.3]{BBR}
A Type II Jordan-Osserman metric is Ricci flat (i.e., $\alpha$ $=$ $\beta$ $=$ $0$)
or otherwise $\alpha$ $=$ $4\beta$ $\neq 0$.
\medskip

\noindent {\normalsize\bf Fact C.} \cite[Proposition 8.4]{BBR}
A Type II Jordan-Osserman metric whose Jacobi operators are not nilpotent
(i.e., $\alpha$ $=$ $4\beta$ $\neq 0$) admits a local parallel field of two-dimensional planes.
\medskip

Therefore, we investigate Walker metrics (i.e., those admitting a locally defined
two-dimensional degenerate parallel distribution) in detail in $\S 2$, with special attention
to the (anti-) self-dual Weyl curvature tensors.  A complete description of self-dual
Walker metrics is given in $\S 2.3$. The integration of the Einstein equation
for a self-dual Walker metric, which lets us determine all pointwise Osserman self-dual Walker metrics,
is carried out in $\S 3$. This leads to (\ref{typeii-1})--(\ref{typeii-2}), thus proving the Main Theorem.

\section{Four-dimensional Walker metrics}\label{Osserman metrics on Walker manifolds}

A \emph{Walker manifold} is a triple $(M,g,\mathcal{D})$ where $M$ is an $n$-dimensional manifold, $g$ an indefinite
metric  and $\mathcal{D}$ an $r$-dimensional parallel null distribution. Of special interest are those
manifolds admitting a field of null planes of maximum dimensionality ($r=\frac{n}{2}$). Since the dimension
of a null plane is $r\leq \frac{n}{2}$, the lowest possible case of a Walker metric is that of
$(++--)$-manifolds admitting  a field of parallel null two-planes.

Observe that there is a tight connection between Walker structures and Osserman metrics. First of all, it is a matter
of fact that all known examples of Osserman metrics with nondiagonalizable Jacobi operators (i.e., Types II
and III) are realized as Walker metrics.
On the other hand,
Walker metrics appear as the underlying structure of several specific pseudo-Riemannian structures.
For instance, the metric tensor of any para-K{\"a}hler \cite{CFG} (and hence any hypersymplectic \cite{Hit}) structure is necessarily
of Walker type. The same occurs for the underlying metric of real hypersurfaces in indefinite
space forms whose shape operator is nilpotent \cite{Ma}.

\subsection{The canonical form of a Walker metric}\label{The canonical form of a Walker metric}

For our purposes it is convenient to use a special coordinate system associated to any Walker metric.
So, let $g$ be a four-dimensional pseudo-Riemannian metric admitting a two-dimensional parallel null
distribution. A canonical form for such a metric has been obtained by Walker \cite{W1} showing the
existence of suitable coordinates $(x_1,\dots,x_4)$ where the metric expresses as
\begin{equation}\label{Walker metric}
g_{(x_1,x_2,x_3,x_4)} =\left(\begin{array}{cccc}
    0&0&1&0\\
    0&0&0&1\\
    1&0&a&c\\
    0&1&c&b\\
\end{array}\right)
\end{equation}
for some functions $a$, $b$ and $c$ depending on the coordinates $(x_1,\dots,x_4)$.
As a matter of notation, throughout this work   we denote by $\partial_i$ the coordinate vectors,
$i=1,\dots,4$. Also, $h_{i_1 \dots i_r}$  means partial derivatives $ \frac{\partial h}{\partial x_{i_1}
\dots\partial x_{i_r}}$, for any function $h(x_1,\dots,x_4)$. Now, a straightforward calculation from
(\ref{Walker metric}) shows that the Levi Civita connection is given by
\begin{equation}\label{Walker metric: Levi-Civita connection}
\,\,
\begin{array}{l}
   \nabla_{\partial_1} \partial_3 = \frac{1}{2} a_1 \partial_1 + \frac{1}{2} c_1 \partial_2,
   \qquad
   \nabla_{\partial_1} \partial_4 = \frac{1}{2} c_1 \partial_1 + \frac{1}{2} b_1 \partial_2,

   \\[0.1in]

   \nabla_{\partial_2} \partial_3 = \frac{1}{2} a_2 \partial_1 + \frac{1}{2} c_2 \partial_2,
   \qquad
   \nabla_{\partial_2} \partial_4 = \frac{1}{2} c_2 \partial_1 + \frac{1}{2} b_2 \partial_2,

   \\[0.1in]

   \nabla_{\partial_3} \partial_3 = \frac{1}{2} ( a a_1 + c a_2 + a_3   )\partial_1 +
                                    \frac{1}{2} ( c a_1 + b a_2  - a_4 + 2 c_3  )\partial_2
   - \frac{a_1}{2}  \partial_3 - \frac{a_2}{2} \partial_4,

   \\[0.1in]

   \nabla_{\partial_3} \partial_4 = \frac{1}{2} ( a_4 + a c_1 + c c_2  )\partial_1
   + \frac{1}{2} ( b_3 + c c_1 + b c_2  )\partial_2
   - \frac{c_1}{2}  \partial_3 - \frac{c_2}{2}  \partial_4,

   \\[0.1in]

   \nabla_{\partial_4} \partial_4 = \frac{1}{2} ( a b_1 + c b_2 - b_3  + 2 c_4  )\partial_1
   + \frac{1}{2} (  c b_1 + b b_2 + b_4   )\partial_2
   -\frac{b_1}{2}  \partial_3 - \frac{ b_2}{2} \partial_4,
\end{array}
\end{equation}
and the Riemann curvature tensor $R(X,Y)=\nabla_{[X,Y]}-[\nabla_X,\nabla_Y]$
satisfies
\begin{equation}\label{Walker metric: curvature tensor R}
\begin{array}{l}
    R_{1313} = -\frac{1}{2} a_{11},\quad
    R_{1314} = -\frac{1}{2} c_{11},\quad
    R_{1323} = -\frac{1}{2} a_{12},\quad
    R_{1324} = -\frac{1}{2} c_{12},
    \\[0.1in]
    R_{1334} = \frac{1}{4} \left(
        -a_2 b_1 + c_1 c_2 + 2 a_{14} - 2 c_{13} \right),

    \\[0.1in]

    R_{1414} = -\frac{1}{2} b_{11},\quad
    R_{1423} = -\frac{1}{2} c_{12},\quad
    R_{1424} = -\frac{1}{2} b_{12},\quad
    \\[0.1in]
    R_{1434} = \frac{1}{4} \left(
        - c_1^2 + a_1 b_1 -b_1 c_2  + b_2 c_1  - 2 b_{13}+ 2 c_{14}
        \right),

    \\[0.1in]

    R_{2323} = -\frac{1}{2} a_{22},\quad
    R_{2324} = -\frac{1}{2} c_{22},\quad
    \\[0.1in]
    R_{2334} = \frac{1}{4} \left(
        c_2^2- a_2 b_2- a_1 c_2  + a_2 c_1  + 2a_{24} - 2c_{23}
        \right),

    \\[0.1in]

    R_{2424} = -\frac{1}{2} b_{22},\quad
    R_{2434} = \frac{1}{4} \left(
            a_2 b_1 - c_1 c_2 - 2 b_{23}+ 2 c_{24}
            \right),

    \\[0.1in]

    R_{3434} = \frac{1}{4} \left(
        -a c_1^2 - b c_2^2
        + a a_1 b_1 + c a_1 b_2- a_1 b_3 +2 a_1 c_4
        \right.
        \\[0.1in]
        \phantom{R_{3434} = \frac{1}{4} \{}
        \left.
        + c a_2 b_1 + b a_2 b_2 + a_2 b_4
        + a_3 b_1
        - a_4 b_2 - 2 a_4 c_1
        \right.
        \\[0.1in]
        \phantom{R_{3434} = \frac{1}{4} \{}
        \left.
        + 2 b_2 c_3
        - 2 b_3 c_2
        - 2 c c_1 c_2
        -2 a_{44} - 2 b_{33} + 4 c_{34}
        \right).
\end{array}
\end{equation}

Next, let $\rho(X,Y)$ $=$ trace$\,\{U\rightsquigarrow R(X,U)Y\}$ and $\tau= \mbox{trace}\,\rho$ be the
Ricci tensor and the scalar curvature, respectively. Then
\begin{equation}\label{Walker Ricci tensor}
\hspace*{0.3cm}
\begin{array}{l}
    \rho_{13} = \frac{1}{2} \left( a_{11} + c_{12} \right),\quad
    \rho_{14} = \frac{1}{2} \left( b_{12} + c_{11} \right),

    \\[0.1in]

    \rho_{23} = \frac{1}{2} \left( a_{12} + c_{22} \right),\quad
    \rho_{24} = \frac{1}{2} \left( b_{22} + c_{12} \right),

    \\[0.1in]

    \rho_{33} = \frac{1}{2} \left(
        -c_2^2
        + a_1   c_2 + a_2  b_2 - a_2  c_1
        + a   a_{11} + 2  c  a_{12}
        + b   a_{22} + 2   c_{23} - 2  a_{24}
        \right),

    \\[0.1in]

    \rho_{34} = \frac{1}{2} \left(
        - a_2  b_1 + c_1   c_2
        + a_{14} + b_{23}
        + a  c_{11} + 2  c  c_{12} - c_{13}
        + b  c_{22} -c_{24}
        \right),

    \\[0.1in]

    \rho_{44} = \frac{1}{2} \left(
        - c_1^2 +
        a_1   b_1 - b_1   c_2 +  b_2   c_1
        + a   b_{11} + 2  c    b_{12} - 2   b_{13} +
        b   b_{22} + 2   c_{14}
        \right),
\end{array}
\end{equation}
and
\begin{equation}\label{curvatura escalar}
    \tau=a_{11} + b_{22} + 2  c_{12}.
\end{equation}

Further, let $W$ denote the Weyl conformal curvature tensor given by
\begin{equation}\label{Weyl curvature tensor}
\begin{array}{rcl}
    W(X,Y,Z,T) & = & R(X,Y,Z,T)
    \\[0.1in]
    && + \frac{\tau}{(n-1)(n-2)} \{g(X,Z)g(Y,T)-g(Y,Z)g(X,T)\}
    \\[0.1in]
    && -\frac{1}{n-2} \{\rho(X,Z)g(Y,T)-\rho(Y,Z)g(X,T)
    \\[0.1in]
    && \hspace*{1.1cm}+ \rho(Y,T)g(X,Z)-\rho(X,T)g(Y,Z)\}.
\end{array}
\end{equation}

\subsection{Self-duality and anti-self-duality conditions}\label{Self-duality and anti-self-duality conditions}

Considering the curvature tensor $R$ as an endomorphism of $\Lambda^2(M)$, we have
the following $O(2,2)$-decomposition
\begin{equation}\label{eq:spliting-1}
R\equiv \frac{\tau}{12}\;Id_{\Lambda^2} +\rho_0 +
W:\,\,\Lambda^2 \rightarrow \Lambda^2
\end{equation}
where $\rho_0$ denotes the traceless Ricci tensor,
$\rho_0(X,Y)$ $=$ $\rho(X,Y)$ $-$ $\frac{\tau}{4}\,g(X,Y)$.
The Hodge star operator $\ast:\Lambda^2\rightarrow \Lambda^2$ associated to
any $(++--)$-metric induces a further splitting
$\Lambda^2 =\Lambda^2_+\oplus\Lambda^2_-$, where $\Lambda^2_\pm$ denotes the $\pm 1$-eigenspaces of the Hodge star operator, that is,
$\Lambda^2_\pm=\{\alpha\in\Lambda^2(M)/ \ast\alpha=\pm\alpha\}$. Correspondingly, the curvature tensor further decomposes as
\begin{equation}\label{eq:spliting-2}
R\equiv \frac{\tau}{12}\;Id_{\Lambda^2} + \rho_0 +
W^++W^-,
\end{equation}
where $W^\pm=\frac{1}{2}\;(W\pm\ast W)$. Recall that a pseudo-Riemannian $4$-manifold is called
\emph{self-dual} (resp., \emph{anti-self-dual}) if $W^-=0$ (resp., $W^+=0$).
The connection between Osserman and (anti-) self-dual metrics relies on Fact A, and thus
 the analysis of the (anti-) self-duality conditions will play a basic role in
what follows.

Let $\{e_1,e_2,e_3,e_4\}$ be an orthonormal basis and, as a convention, assume that $e_1$ and $e_2$ are
spacelike while $e_3$ and $e_4$ are timelike vectors.
Now, local bases of the spaces of self-dual and anti-self-dual two-forms can be constructed as
$$
\Lambda^2_\pm =\langle\left\{ E_1^\pm, E_2^\pm, E_3^\pm\right\}\rangle,
$$
where
$$
E_1^\pm = \frac{e^1\wedge e^2 \pm e^3\wedge e^4}{\sqrt{2}},\quad
E_2^\pm = \frac{e^1\wedge e^3 \pm e^2\wedge e^4}{\sqrt{2}},\quad
E_3^\pm = \frac{e^1\wedge e^4 \mp e^2\wedge e^3}{\sqrt{2}}.
$$
Here observe that the Hodge star operator satisfies
$$
e^i\wedge e^j\wedge \star(e^k\wedge e^l) =
    (\delta^i_k \delta^j_l-\delta^i_l \delta^j_k)\:\varepsilon_{i}\varepsilon_{j}\:
    e^1 \wedge e^2\wedge e^3\wedge e^4,
$$
where $\varepsilon_{i}=g(e_i,e_i)$. Further note that
$\langle E_1^\pm, E_1^\pm \rangle =1$, $\langle E_2^\pm, E_2^\pm \rangle =-1$,
$\langle E_3^\pm, E_3^\pm \rangle =-1$, and therefore  with respect to the above bases the
self-dual and anti-self-dual Weyl curvature operators
$W^\pm:\Lambda^2_\pm \longrightarrow \Lambda^2_\pm$ have the following matrix form:
\begin{equation}\label{W+- matrix form}
    W^\pm =
    \left(
    \begin{array}{rrr}
         W_{11}^\pm &
         W_{12}^\pm &
         W_{13}^\pm
        \\[0.1in]
        - W_{12}^\pm&
        - W_{22}^\pm&
        - W_{23}^\pm
        \\[0.1in]
        - W_{13}^\pm &
        - W_{23}^\pm &
        - W_{33}^\pm
    \end{array}
    \right),
\end{equation}
where $W^\pm_{ij}$ $=$ $W(E_i^\pm, E_j^\pm)$ and $W(e^i\wedge e^j, e^k\wedge e^l)$ $=$ $W(e_i,e_j,e_k,e_l)$.

Next, for a Walker metric (\ref{Walker metric}) an orthonormal basis can be specialized by
using the canonical coordinates as follows:
\begin{equation}\label{orthonormal basis}
\begin{array}{ll}
   e_1 = \frac{1}{2} (1-a) \partial_1 + \partial_3,  &
   e_2 = - c \partial_1 + \frac{1}{2} (1-b) \partial_2 + \partial_4,

   \\[0.1in]

   e_3 = -\frac{1}{2} (1+a) \partial_1 + \partial_3,  \qquad &
   e_4 = - c \partial_1 - \frac{1}{2} (1+b) \partial_2 + \partial_4.
\end{array}
\end{equation}
Now, a long but straightforward calculation using (\ref{Walker metric: curvature tensor R}), (\ref{Walker Ricci
tensor}) and (\ref{curvatura escalar}) shows that the components of $W^-$ in (\ref{W+- matrix form}) are given by
\begin{equation}\label{components of W-}
\begin{array}{ll}
   W^-_{11} = -\frac{1}{12} ( a_{11} + 3 a_{22} + 3 b_{11} + b_{22} - 4 c_{12} ),
   \\[0.1in]
   W^-_{22} = -\frac{1}{6} ( a_{11} + b_{22} - 4 c_{12} ),
   \\[0.1in]
   W^-_{33} = \frac{1}{12} ( a_{11} - 3 a_{22} - 3 b_{11} + b_{22} - 4 c_{12} ),
   \\[0.1in]
   W^-_{12} = \frac{1}{4}( a_{12} + b_{12} - c_{11} - c_{22} ),
   \\[0.1in]
   W^-_{13} = \frac{1}{4} ( a_{22} - b_{11} ),
   \\[0.1in]
   W^-_{23} = -\frac{1}{4} ( a_{12} - b_{12} + c_{11} - c_{22} ),
\end{array}
\end{equation}
while the components of $W^+$ are determined by $W^+_{11}$, $W^+_{12}$ and the scalar curvature as follows:
\begin{equation}\label{components of W+}
\qquad
   W^+_{22} = -\frac{\tau}{6},  \qquad
   W^+_{33} = W^+_{11} + \frac{\tau}{6},  \qquad
   W^+_{13} = W^+_{11} + \frac{\tau}{12}, \qquad
   W^+_{23} = W^+_{12},
\end{equation}
the expressions of $W^+_{11}$ and $W^+_{12}$ being
\begin{equation}\label{W+ 11}
\,
\begin{array}{l}
   W^+_{11} = \frac{1}{12} (
   6 c a_1 b_2 - 6 a_1 b_3 - 6 b a_1 c_2  + 12 a_1 c_4 - 6 c a_2 b_1 + 6 a_2 b_4 + 6 b a_2 c_1
   \\[0.1in]
   \phantom{W^+_{11} = \frac{1}{12} ( }
   +  6 a_3 b_1 - 6 a_4 b_2 - 12 a_4 c_1  + 6 a b_1 c_2 - 6 a b_2 c_1 + 12 b_2 c_3 - 12 b_3 c_2

   \\[0.1in]
   \phantom{W^+_{11} = \frac{1}{12} ( }
   - a_{11} - 12 c^2 a_{11} - 12 b c a_{12} + 24 c a_{14} - 3 b^2 a_{22} + 12 b a_{24} -12 a_{44}
   \\[0.1in]
   \phantom{W^+_{11} = \frac{1}{12} ( }
   - 3 a^2 b_{11} + 12 a b_{13} - b_{22} - 12 b_{33} + 12 a c c_{11} - 2 c_{12} + 6 a b c_{12}

   \\[0.1in]
   \phantom{W^+_{11} = \frac{1}{12} ( }
   - 24 c c_{13} - 12 a c_{14} - 12 b c_{23}  + 24 c_{34}
   ),
\end{array}
\end{equation}
\medskip
\begin{equation}\label{W+ 12}
\begin{array}{l}
   W^+_{12} = \frac{1}{4}  (
   - 2 c a_{11} - b a_{12} + 2 a_{14} + a b_{12} - 2 b_{23} + a c_{11} - 2 c c_{12} - 2 c_{13}
   \\[0.1in]
   \phantom{W^+_{12} = \frac{1}{4}  (  }
   - b c_{22} + 2 c_{24}
   ).
\end{array}
\end{equation}

\begin{remark}\label{re: autovalores de W+ en general}\rm
It is important to recognize here that the connection between Einstein (anti-) self-dual
and pointwise Osserman manifolds at Fact A goes further to the Jordan normal forms of the
nonzero part of the Weyl curvature tensor $W^\pm$
and the Jacobi operators (cf. \cite{GKV}).
So pointwise Osserman manifolds whose Jacobi operators are of Type Ia, Ib, II or III
correspond to self-dual (or anti-self-dual) Einstein manifolds whose self-dual
(or anti-self-dual) Weyl curvature tensor is of Type Ia, Ib, II or III, respectively
(see also \cite{BGi}).

The relations in (\ref{components of W+}) among the components of $W^+$ show that
$$
W^+ = \left(
\begin{array}{ccc}
   W^+_{11} & W^+_{12} & W^+_{11} + \frac{\tau}{12}
   \\[0.1in]
   - W^+_{12} & \frac{\tau}{6} & - W^+_{12}
   \\[0.1in]
   -(W^+_{11} + \frac{\tau}{12}) & - W^+_{12} & -(W^+_{11} + \frac{\tau}{6})
\end{array}
\right),
$$
and, as a consequence, the eigenvalues of $W^+$ are given by
\begin{equation}\label{eigenvalues W^+}
   \left\{  \frac{\tau}{6}, -\frac{\tau}{12}, -\frac{\tau}{12}  \right\}.
\end{equation}
Since the induced metric on $\Lambda^2_+$ is of Lorentzian signature
the structure of $W^+$ is determined by its Jordan normal form,
which may correspond to Type Ia or type II/III, depending on whether $W^+$ is diagonalizable
or not. A straightforward calculation shows that
$$
   \left( W^+ - \frac{\tau}{6} Id \right) \cdot \left( W^+ + \frac{\tau}{12} Id \right) =
   \frac{\tau^2 + 12 \tau W^+_{11} + 48 \left(W^+_{12}\right)^2}{48}
   \left(
   \begin{array}{ccc}
      -1 & 0 & -1
      \\
      0 & 0 & 0
      \\
      1 & 0 & 1
   \end{array}
   \right),
$$
from where we have the following:
\begin{itemize}
\item[(i)] If $\tau\neq 0$, then $W^+$ has  nonzero eigenvalues $\left\{  \frac{\tau}{6}, -\frac{\tau}{12},
-\frac{\tau}{12}  \right\}$, and
\end{itemize}
\begin{equation}\label{diagonalization of W+}
   \tau^2 + 12 \tau W^+_{11} + 48 \left(W^+_{12}\right)^2 = 0
\end{equation}
\begin{itemize}
\item[]
is the necessary and sufficient condition for the diagonalizability of $W^+$. If (\ref{diagonalization of W+})
does not hold, then $-\frac{\tau}{12}$ is a double root of the minimal polynomial of $W^+$.

\item[(ii)] If $\tau=0$, then $W^+$ vanishes if and only if $W^+_{11}=W^+_{12}=0$ and moreover
   \begin{itemize}
      \item[(ii.1)] $W^+$ is two-step nilpotent if and only if $W^+_{11}\neq 0$ and $W^+_{12}=0$,
      \item[(ii.2)] $W^+$ is three-step nilpotent if and only if   $W^+_{12}\neq 0$.
   \end{itemize}
\end{itemize}
On the other hand, observe from (\ref{eigenvalues W^+}) that any anti-self-dual Walker metric has vanishing scalar curvature
and hence Einstein anti-self-dual Walker metrics are Ricci flat.
\end{remark}

\subsection{Explicit form of Self-dual Walker metrics}\label{Section: self-dual Walker metrics}

Recall that our main purpose is to obtain a description of non Ricci flat Type II Jordan-Osserman
four-manifolds. Then, as a consequence of Remark \ref{re: autovalores de W+ en general}, we
restrict our analysis to self-dual Walker metrics. Thus, in what remains of this section
we will give a complete description of self-dual Walker metrics by integrating the PDE system
given by (\ref{components of W-}).

\begin{theorem}\label{th:self-dual Walker metrics}
A Walker metric {\rm (\ref{Walker metric})} is self-dual if and only if the defining functions
$a(x_1,x_2,x_3,x_4)$, $b(x_1,x_2,x_3,x_4)$ and $c(x_1,x_2,x_3,x_4)$ are given by
\begin{equation}\label{a, b and c for self-dual Walker metrics}
\begin{array}{l}
   a(x_1,x_2,x_3,x_4) = x_1^3 \mathcal{A} + x_1^2 \mathcal{B} + x_1^2 x_2 \mathcal{C}
   + x_1 x_2 \mathcal{D} + x_1 P + x_2 Q + \xi,
   \\[0.1in]
   b(x_1,x_2,x_3,x_4) = x_2^3 \mathcal{C} + x_2^2 \mathcal{E} + x_1 x_2^2 \mathcal{A}
   + x_1 x_2 \mathcal{F} + x_1 S + x_2 T + \eta,
   \\[0.1in]
   c(x_1,x_2,x_3,x_4) = \frac{1}{2} x_1^2 \mathcal{F} + \frac{1}{2} x_2^2 \mathcal{D}
   + x_1^2 x_2 \mathcal{A} + x_1 x_2^2 \mathcal{C} + \frac{1}{2} x_1 x_2 ( \mathcal{B} + \mathcal{E} )
   \\[0.1in]
   \phantom{c(x_1,x_2,x_3,x_4) =}
   + x_1 U + x_2 V + \gamma,
\end{array}
\end{equation}
where capital, calligraphic and Greek letters are all arbitrary smooth functions depending only on the
coordinates $(x_3,x_4)$.
\end{theorem}

\noindent{\bf Proof.} First of all observe that using (\ref{components of W-}) the  self-duality can be
initially characterized by means of the following  five PDEs:
\begin{equation}\label{self-duality PDEs}
   \begin{array}{l}
      a_{22} = 0,
      \\[0.1in]
      b_{11} = 0,
      \\[0.1in]
      a_{12} - c_{22} = 0,
      \\[0.1in]
      b_{12} - c_{11} = 0,
      \\[0.1in]
      a_{11} + b_{22} - 4 c_{12} = 0.
   \end{array}
\end{equation}
So (\ref{a, b and c for self-dual Walker metrics}) is obtained as the solution of (\ref{self-duality PDEs}).
We proceed.

\medskip

\noindent\emph{The first step.} First and second equations in (\ref{self-duality PDEs}) translate into
\begin{equation}\label{self-duality 1}
\begin{array}{l}
   a(x_1,x_2,x_3,x_4) = x_2  A(x_1,x_3,x_4) + B(x_1,x_3,x_4),
   \\[0.1in]
   b(x_1,x_2,x_3,x_4) = x_1  C(x_2,x_3,x_4) + D(x_2,x_3,x_4)
\end{array}
\end{equation}
and hence the third equation in (\ref{self-duality PDEs}), i.e., $c_{22}(x_1,x_2,x_3,x_4)=A_1(x_1,x_3,x_4)$,
implies that
\begin{equation}\label{self-duality 2}
   c(x_1,x_2,x_3,x_4) = \frac{1}{2} x_2^2  A_1(x_1,x_3,x_4)
                      + x_2 E(x_1,x_3,x_4) + F(x_1,x_3,x_4).
\end{equation}

\medskip

\noindent\emph{The second step.} In this step we make a long but straightforward use of the fourth equation
in (\ref{self-duality PDEs}), which by (\ref{self-duality 1}) and (\ref{self-duality 2}) means
\begin{equation}\label{self-duality 4}
\begin{array}{l}
   \frac{1}{2} x_2^2 A_{111}(x_1,x_3,x_4) + x_2 E_{11}(x_1,x_3,x_4)
   \\[0.1in]
   \phantom{\frac{1}{2} x_2^2 A_{111}(x_1,x_3,x_4)}
   + F_{11}(x_1,x_3,x_4) -    C_{2}(x_2,x_3,x_4)=0.
\end{array}
\end{equation}
First, appropriate successive differentiations in the above equation lead to $A_{1111}(x_1,$ $x_3,$ $x_4) =
0$, $E_{111}(x_1,x_3,x_4)=0$ and $F_{111}(x_1,x_3,x_4)=0$, and hence
\begin{equation}\label{self-duality 5}
\begin{array}{l}
   A(x_1,x_3,x_4) = x_1^3 G(x_3,x_4) + x_1^2 \mathcal{C}(x_3,x_4) + x_1 \mathcal{D}(x_3,x_4) + Q(x_3,x_4),
   \\[0.1in]
   E(x_1,x_3,x_4) = x_1^2 H(x_3,x_4) + x_1 I(x_3,x_4) + V(x_3,x_4),
   \\[0.1in]
   F(x_1,x_3,x_4)=x_1^2 J(x_3,x_4) + x_1 U(x_3,x_4) + \gamma(x_3,x_4).
\end{array}
\end{equation}
By (\ref{self-duality 5}), condition (\ref{self-duality 4}) reduces to
\begin{equation}\label{self-duality 8}
   3 x_2^2 G(x_3,x_4) + 2 x_2 H(x_3,x_4) + 2 J(x_3,x_4) - C_{2}(x_2,x_3,x_4)=0,
\end{equation}
from where $C_{222}(x_2,x_3,x_4)=6 G(x_3,x_4)$, which implies that
\begin{equation}\label{self-duality 9}
   C(x_2,x_3,x_4) = x_2^3  G(x_3,x_4) + x_2^2 \mathcal{A}(x_3,x_4) + x_2 \mathcal{F}(x_3,x_4) + S(x_3,x_4).
\end{equation}

The final form of (\ref{self-duality 8}) is given by
$$
   2 x_2 ( \mathcal{A}(x_3,x_4) - H(x_3,x_4) )  - 2 J(x_3,x_4) +   \mathcal{F}(x_3,x_4)= 0,
$$
and as a consequence
\begin{equation}\label{self-duality 11}
   H(x_3,x_4) = \mathcal{A}(x_3,x_4), \qquad
   J(x_3,x_4) = \frac{1}{2} \mathcal{F}(x_3,x_4).
\end{equation}
Collecting together the information in (\ref{self-duality 5}), (\ref{self-duality 9}) and (\ref{self-duality
11}), at the end of this step we have that (\ref{self-duality 1}) and (\ref{self-duality 2}) transform into
\begin{equation}\label{self-duatilty step 2}
\begin{array}{l}
   a(x_1,x_2,x_3,x_4) = x_1^3 x_2 G + x_1^2 x_2 \mathcal{C} + x_1 x_2 \mathcal{D} + x_2 Q + B(x_1,x_3,x_4),
   \\[0.1in]
   b(x_1,x_2,x_3,x_4) = x_1 x_2^3 G + x_1 x_2^2 \mathcal{A} + x_1 x_2 \mathcal{F} + x_1 S + D(x_2,x_3,x_4),
   \\[0.1in]
   c(x_1,x_2,x_3,x_4) = \frac{3}{2} x_1^2 x_2^2 G + \frac{1}{2} x_1^2 \mathcal{F} + \frac{1}{2} x_2^2 \mathcal{D}
   + x_1^2 x_2 \mathcal{A} + x_1 x_2^2 \mathcal{C}
   \\[0.1in]
   \phantom{c(x_1,x_2,x_3,x_4) = }
   + x_1 x_2 I + x_1 U + x_2 V + \gamma.
\end{array}
\end{equation}

\medskip

\noindent\emph{The third step.} To finish the determination of the defining functions we deal with the last
equation in (\ref{self-duality PDEs}), which under the expressions in (\ref{self-duatilty step 2}) and
differentiating by $x_1$ and $x_2$ leads to
\begin{equation}\label{self-duality 13}
   G(x_3,x_4)=0,
\end{equation}
and hence that equation reduces to
\begin{equation}\label{self-duality 14}
\begin{array}{l}
   6 x_1 \mathcal{A}(x_3,x_4)+ 6 x_2 \mathcal{C}(x_3,x_4) + 4 I(x_3,x_4)
   \\[0.1in]
   \phantom{6 x_1 \mathcal{A}(x_3,x_4)}
   - B_{11}(x_1,x_3,x_4) - D_{22}(x_2,x_3,x_4) = 0.
\end{array}
\end{equation}
It follows that $B_{1111}(x_1,x_3,x_4)=0$, from where
\begin{equation}\label{self-duality 15}
   B(x_1,x_3,x_4) = x_1^3  L(x_3,x_4) + x_1^2 \mathcal{B}(x_3,x_4) + x_1 P(x_3,x_4) + \xi(x_3,x_4)
\end{equation}
and (\ref{self-duality 14}) takes the form
\begin{equation}\label{self-duality 16}
\begin{array}{l}
   6 x_1 (\mathcal{A}(x_3,x_4) - L(x_3,x_4)) + 6 x_2 \mathcal{C}(x_3,x_4)
   \\[0.1in]
   \phantom{6 x_1 \mathcal{A}(x_3,x_4)}
   + 4 I(x_3,x_4) - 2 \mathcal{B}(x_3,x_4)  - D_{22}(x_2,x_3,x_4) = 0,
\end{array}
\end{equation}
which leads to
\begin{equation}\label{self-duality 17}
   L(x_3,x_4)=\mathcal{A}(x_3,x_4).\qquad
\end{equation}
Thus, by (\ref{self-duality 13}), (\ref{self-duality 15}) and (\ref{self-duality 17}), the expression of $a$
in (\ref{self-duatilty step 2}) transforms  into
\begin{equation}\label{self-duality final a}
   a(x_1,x_2,x_3,x_4) = x_1^3 \mathcal{A} + x_1^2 \mathcal{B} + x_1^2 x_2 \mathcal{C}
   + x_1 x_2 \mathcal{D} + x_1 P + x_2 Q + \xi,
\end{equation}
which finishes the process for this defining function. At this point, (\ref{self-duality 16}) has  the form
\begin{equation}\label{self-duality 18}
   6 x_2 \mathcal{C}(x_3,x_4)  + 4 I(x_3,x_4) - 2 \mathcal{B}(x_3,x_4)  - D_{22}(x_2,x_3,x_4) = 0,
\end{equation}
from where $D_{2222}(x_2,x_3,x_4)=0$, and hence
\begin{equation}\label{self-duality 19}
   D(x_2,x_3,x_4) = x_2^3  M(x_3,x_4) + x_2^2 \mathcal{E}(x_3,x_4) + x_2 T(x_3,x_4) + \eta(x_3,x_4).
\end{equation}
Then, (\ref{self-duality 18}) reduces to
$$
   3 x_2 ( \mathcal{C}(x_3,x_4)-M(x_3,x_4) )   + 2 I(x_3,x_4) -  \mathcal{B}(x_3,x_4)  - \mathcal{E}(x_3,x_4) = 0,
$$
so we conclude that
\begin{equation}\label{self-duality 21}
   M(x_3,x_4)=\mathcal{C}(x_3,x_4),\qquad
   I(x_3,x_4)=\frac{1}{2} (\mathcal{B} + \mathcal{E}).
\end{equation}
Thus, (\ref{self-duality 13}), (\ref{self-duality 19}) and (\ref{self-duality 21}) show  that $b$ and $c$ in
(\ref{self-duatilty step 2}) take the desired form
\begin{equation}\label{self-duality final b}
\begin{array}{l}
   b(x_1,x_2,x_3,x_4) = x_2^3 \mathcal{C} + x_2^2 \mathcal{E} + x_1 x_2^2 \mathcal{A}
   + x_1 x_2 \mathcal{F} + x_1 S + x_2 T + \eta,
   \\[0.1in]
   c(x_1,x_2,x_3,x_4) = \frac{1}{2} x_1^2 \mathcal{F} + \frac{1}{2} x_2^2 \mathcal{D}
   + x_1^2 x_2 \mathcal{A} + x_1 x_2^2 \mathcal{C} + \frac{1}{2} x_1 x_2 ( \mathcal{B} + \mathcal{E} )
   \\[0.1in]
   \phantom{c(x_1,x_2,x_3,x_4) =}
   + x_1 U + x_2 V + \gamma,
\end{array}
\end{equation}
finishing the proof.$\hfill\square$

\begin{remark}\rm
A four-dimensional Walker metric is said to be \emph{strict} if it admits two orthogonal parallel null
vector fields rather than a parallel two-dimensional null distribution. It follows from the work of Walker
that any strict Walker metric is given by (\ref{Walker metric}) for any functions $a$, $b$ and $c$
depending only on the coordinates $(x_3,x_4)$. Now  (\ref{Walker Ricci tensor}) and
Theorem \ref{th:self-dual Walker metrics} imply that any strict Walker metric is Ricci flat and self-dual, and
hence Osserman. Moreover, Remark \ref{re: autovalores de W+ en general} shows that the Jacobi operators
are identically zero or otherwise they are two-step nilpotent (depending on whether $W_{11}^+=2c_{34}-a_{44}-b_{33}$
vanishes or not).
\end{remark}

\section{Proof of the Main Theorem}

Recall that our purpose is to obtain a local description of Type II Jordan-Osserman metrics whose
Jacobi operators have nonzero eigenvalues. In such a case, the eigenvalues must be in a ratio
$1$ : $\frac{1}{4}$ and the underlying metric is a Walker metric (cf. Facts B and C).
Therefore in order to achieve the desired result,
only Osserman metrics on Walker manifolds deserve further consideration and moreover, it immediately
follows from Remark \ref{re: autovalores de W+ en general} that we may restrict to those being self-dual.
In what follows we obtain a complete description of self-dual Einstein Walker metrics in Theorem \ref{th:Einstein self-dual Walker metrics},
from where our main result is derived.

\begin{theorem}\label{th:Einstein self-dual Walker metrics}
A Walker metric {\rm (\ref{Walker metric})} is pointwise Osserman self-dual if and only if
one of the following holds:
\begin{enumerate}
\item[(i)]
The scalar curvature $\tau$ is nonzero and the metric tensor is completely determined by
the functions $a(x_1,\dots, x_4)$, $b(x_1,\dots, x_4)$ and $c(x_1,\dots, x_4)$ as follows
\end{enumerate}
\begin{equation}\label{Einstein self-dual Walker metrics - Sc nonnull: a,b,c}
\begin{array}{l}
   a = x_1^2 \frac{\tau}{6} + x_1 P + x_2 Q
   + \frac{6}{\tau} \left\{ Q   ( T - U ) + V   (P - V) - 2 (Q_4 -  V_3) \right\},
   \\[0.1in]
   b = x_2^2 \frac{\tau}{6} + x_1 S + x_2 T
   + \frac{6}{\tau} \left\{ S   ( P - V ) + U   (T - U) - 2 (S_3 - U_4) \right\},
   \\[0.1in]
   c = x_1 x_2 \frac{\tau}{6} + x_1 U + x_2 V
   + \frac{6}{\tau} \left\{ - Q S + U V + T_3 - U_3 + P_4 - V_4 \right\},
\end{array}
\end{equation}
\begin{enumerate}
\item[]
where capital letters are arbitrary functions depending only on $(x_3, x_4)$.

The Jacobi operators have eigenvalues $\{0, \frac{\tau}{6}, \frac{\tau}{24}, \frac{\tau}{24}\}$ and
they are diagonalizable if and only if {\rm (\ref{diagonalization of W+})} holds. Otherwise,
$\frac{\tau}{24}$ is a double root of the minimal polynomial of the Jacobi operators and the Walker manifold
is Jordan-Osserman on the open set where {\rm (\ref{diagonalization of W+})} does not hold.

\item[(ii)]
The scalar curvature vanishes and the metric tensor is given by
\end{enumerate}
\begin{equation}\label{Einstein self-dual Walker metrics - Sc=0: a,b,c}
\begin{array}{l}
   a =   x_1 P + x_2 Q + \xi,
   \\[0.1in]
   b  =   x_1 S + x_2 T + \eta,
   \\[0.1in]
   c =   x_1 U + x_2 V + \gamma,
\end{array}
\end{equation}
\begin{enumerate}
\item[]
where $P$, $Q$, $S$, $T$, $U$, $V$, $\xi$, $\eta$ and $\gamma$ are smooth functions depending only on $(x_3,x_4)$ satisfying
\end{enumerate}
\begin{equation}\label{Einstein self-dual PDEs}
\begin{array}{l}
      2 (Q_4 -  V_3) = Q   ( T - U ) + V   (P - V),
   \\[0.1in]
      2 (S_3 - U_4) =  S   ( P - V ) + U   (T - U) ,
   \\[0.1in]
     T_3 - U_3 + P_4 - V_4 = Q S - U V.
\end{array}
\end{equation}
\begin{enumerate}
\item[]
In this case, the Jacobi operators have zero eigenvalues and their minimal
polynomials satisfy:
\begin{itemize}
  \item[(ii.1)] The Jacobi operators are vanishing (i.e., Type Ia) if and only if
\end{itemize}
\end{enumerate}
  \begin{equation}\label{W+ 12  si Sc=0}
     T_3 + U_3 - P_4 - V_4 = 0
  \end{equation}
  \begin{enumerate}
  \item[]
\begin{itemize}
\item[]
  and $W^+_{11}$ given by {\rm (\ref{W+ 11})} is also null.

  \item[(ii.2)] The Jacobi operators are two-step nilpotent (i.e., Type II) if and only if {\rm (\ref{W+ 12  si Sc=0})} holds and
  $W^+_{11}$ given by {\rm (\ref{W+ 11})} is nonnull.

  \item[(ii.3)] The Jacobi operators are three-step nilpotent (i.e. Type III) if and only if {\rm (\ref{W+ 12  si Sc=0})} does not
  hold.
\end{itemize}
\end{enumerate}
\end{theorem}

\noindent{\bf Proof.}
Let $\rho_0$ $=$ $\rho-\frac{\tau}{4}g$ be the trace-free Ricci tensor. Then the Einstein equations for a
general Walker metric (\ref{Walker metric}) are as follows.

\begin{equation}\label{Walker Einstein PDEs}
   \begin{array}{l}
      (\rho_0)_{13}=-(\rho_0)_{24}=(\rho_0)_{31}=-(\rho_0)_{42}=\frac{1}{4} \left( a_{11}-b_{22} \right) = 0,
      \\[0.1in]
      (\rho_0)_{14}=(\rho_0)_{41}=\frac{1}{2} \left( b_{12}+c_{11} \right) = 0,
      \\[0.1in]
      (\rho_0)_{23}=(\rho_0)_{32}=\frac{1}{2} \left( a_{12}+c_{22} \right) = 0,
      \\[0.1in]

      (\rho_0)_{33}=\frac{1}{4} \left(  2 a_1 c_2 + 2 a_2 b_2 - 2 a_2 c_1 -2 c_2^2 + a a_{11}     \right.
      \\[0.1in]
      \phantom{(\rho_0)_{33}=\frac{1}{4} (}
      \left. + 4 c a_{12} + 2 b a_{22} - 4a_{24}-a b_{22} - 2 a c_{12} + 4c_{23}   \right)=0,

      \\[0.1in]
      (\rho_0)_{34}=(\rho_0)_{43}=\frac{1}{4}
      \left( - 2 a_2 b_1 + 2 c_1 c_2 - c a_{11} + 2 a_{14} - c b_{22}
        \right.
      \\[0.1in]
      \phantom{(\rho_0)_{34}=(\rho_0)_{43}=\frac{1}{4} (}
      \left.  + 2 b_{23} + 2 a c_{11}+ 2 c c_{12} - 2 c_{13} + 2 b c_{22}
      -2c_{24}       \right)=0,

      \\[0.1in]
      (\rho_0)_{44}=\frac{1}{4}
      \left( 2 a_1 b_1 - 2  b_1 c_2 + 2 b_2 c_1  - 2 c_1^2 - b a_{11}    \right.
      \\[0.1in]
      \phantom{(\rho_0)_{34} =\frac{1}{4} (}
      \left. + 2 a b_{11} + 4 c b_{12} - 4 b_{13}
       + b b_{22}- 2 b c_{12} + 4 c_{14}    \right)=0.
   \end{array}
\end{equation}
Now, since the manifold is self-dual, the defining functions are completely determined by
(\ref{a, b and c for self-dual Walker metrics}) (cf. Theorem \ref{th:self-dual Walker metrics}). Then,
computing the first three equations in (\ref{Walker Einstein PDEs}), we get
$$
     2 x_1 \mathcal{A} - 2 x_2 \mathcal{C} + \mathcal{B} -
   \mathcal{E}=0,\qquad
     2 x_2 \mathcal{A} + \mathcal{F}=0,\qquad
    2 x_1 \mathcal{C} + \mathcal{D}=0,
$$
from where
\begin{equation}\label{Einstein self-duality 1}
   \mathcal{A} = \mathcal{C} = \mathcal{D} = \mathcal{F} = 0, \qquad
   \mathcal{E} = \mathcal{B}.
\end{equation}
With these conditions, and using (\ref{curvatura escalar}), the (constant) scalar
curvature is given by
$$
   \tau = 6 \mathcal{B},
$$
and hence (\ref{a, b and c for self-dual Walker metrics}) reduces to
\begin{equation}\label{Einstein self-duality 2}
\begin{array}{l}
   a  =  x_1^2 \frac{\tau}{6} + x_1 P + x_2 Q + \xi,
   \\[0.1in]
   b  =  x_2^2 \frac{\tau}{6}  + x_1 S + x_2 T + \eta,
   \\[0.1in]
   c  =   x_1 x_2 \frac{\tau}{6} + x_1 U + x_2 V + \gamma.
\end{array}
\end{equation}
Now the two cases are obtained just observing that the last three equations in (\ref{Walker Einstein PDEs})
transform into
$$
\begin{array}{l}
    \frac{\tau}{6} \xi -  \left\{ Q   ( T - U ) + V   (P - V) - 2 (Q_4 -  V_3)
   \right\}=0,

   \\[0.1in]

    \frac{\tau}{6} \gamma -  \left\{ - Q S + U V + T_3 - U_3 + P_4 - V_4 \right\}=0,

    \\[0.1in]

    \frac{\tau}{6} \eta -   \left\{ S   ( P - V ) + U   (T - U) - 2 (S_3 - U_4) \right\}=0,
\end{array}
$$
and noting that when the scalar curvature $\tau$ does not vanish we can determine $\xi$, $\eta$ and $\gamma$
from above, while if $\tau=0$ then we get (\ref{Einstein self-dual PDEs}).

Finally, the eigenvalues and the minimal polynomial of the Jacobi operators for the two cases are obtained as
a direct application of Remark \ref{re: autovalores de W+ en general}, since the eigenvalues and the minimal
polynomial of the self-dual Weyl tensor $W^+$ determine the behavior of the Jacobi operators of a pointwise Osserman
self-dual manifold.$\hfill\square$

\medskip

As a consequence of the previous theorem  and
Remark \ref{re: autovalores de W+ en general}, we have the following characterization of Jordan-Osserman Walker metrics.

\begin{theorem}
Let $(M,g)$ be a Jordan-Osserman Walker $4$-manifold. Then one of the following holds:
\begin{enumerate}
\item[(i)] If the Jacobi operators are diagonalizable, then
$(M,g)$ is either flat or locally isometric to a paracomplex space form.
\item[(ii)] If the Jacobi operators are non diagonalizable, then either
\begin{enumerate}
\item[(ii.1)] the Jacobi operators are two-step or three-step nilpotent
\item[(ii.2)] the metric $g$ is given by {\rm (\ref{typeii-1})--(\ref{typeii-2})}.
\end{enumerate}
\end{enumerate}
\end{theorem}

\noindent{\bf Proof.}
Four-dimensional Jordan-Osserman manifolds with diagonalizable Jacobi operators have been classified at
\cite{BBR}, showing that they correspond to real, complex or paracomplex space forms. Next note that real and complex space forms
do not support a Walker metric unless they are flat. Indeed, any space
of constant curvature is locally conformally flat and hence vanishing of $W^+$ shows that any such a Walker metric must be flat.
Analogously, K{\"a}hler metrics of constant holomorphic
sectional curvature have zero Bochner tensor, which shows that $W^+=0$ \cite{Br}, \cite{GVR}.
Hence, no K{\"a}hler metric of constant holomorphic sectional curvature  may be Walker unless it is flat.

On the other hand, if the Jacobi operators are nondiagonalizable, they must be of Type II or III, since Type Ib cannot occur \cite{BBR}.
Then, since anti-self-dual Jordan-Osserman Walker metrics have vanishing scalar curvature, the corresponding Jacobi operators are
either two-step or three-step nilpotent. The only remaining case is that of self-dual Jordan-Osserman Walker metrics, which corresponds
to the Main Theorem, thus finishing the proof.$\hfill\square$

\begin{remark}\rm
Para-K{\"a}hler manifolds of constant paraholomorphic sectional curvature $\alpha$ can  be easily
described as Walker manifolds. For instance, let $a$, $b$ and $c$ be the coordinate functions
given by
$$
a(x_1,x_2,x_3,x_4)=\alpha x_1^2 , \quad b(x_1,x_2,x_3,x_4)=\alpha x_2^2 , \quad c(x_1,x_2,x_3,x_4)=\alpha x_1x_2,
$$
and $J$ the paracomplex structure
$$
J\partial_1 =-\partial_1,\quad J\partial_2 =-\partial_2, \quad J\partial_3=-a\partial_1-c\partial_2+\partial_3,
\quad J\partial_4=-c\partial_1-b\partial_2+\partial_4.
$$
Then $(\mathbb{R}^4,g,J)$ is a para-K{\"a}hler manifold of constant paraholomorphic sectional curvature $\alpha$.
\end{remark}

\begin{remark}\rm
Notice from Remark \ref{re: autovalores de W+ en general} that any  anti-self-dual Jordan-Osserman Walker metric has necessarily
nilpotent Jacobi operators. However besides the fact that many nilpotent Jordan-Osserman metrics are known, none
of the previous examples are anti-self-dual but all of them correspond to special cases of
Theorem \ref{th:Einstein self-dual Walker metrics}.
The general expression of $W^+_{11}$ at (\ref{W+ 11}) makes it quite untractable and hence it is very difficult to obtain the general
form of anti-self-dual Walker metrics. However, for the special choice of
$a(x_1,x_2,x_3,x_4)$ $=$ $b(x_1,x_2,x_3,x_4)$ $=$ $c(x_1,x_2,x_3,x_4)$, anti-self-dual Einstein
metrics are characterized by
$$
\begin{array}{l}
a_{11}=a_{22}=-a_{12},\\
a_{13}=a_{14},\\
a_{23}=a_{24},\\
a_{33}+a_{44}=2a_{34}.
\end{array}
$$
Now, it follows that $a(x_1,x_2,x_3,x_4)$ $=$ $x_1^2+x_2^2-2x_1x_2$ defines an Osserman anti-self-dual Walker metric
with two-step nilpotent Jacobi operators.
\end{remark}

\begin{remark}\rm
Note that any Type III Jordan-Osserman Walker metric is Ricci flat, and thus
the Jacobi operators are three-step nilpotent. The existence of non-Ricci flat Type III
metrics is still an open problem.
\end{remark}

\begin{remark}\rm
As a final remark, recall that when
dealing with Osserman metrics, main attention is usually paid to the behavior of the
eigenvalues of the Jacobi operators. However when the metric under consideration is of
indefinite signature, the corresponding eigenspaces play a basic role too. Indeed,
four-dimensional complex and
paracomplex space forms have diagonalizable Jacobi operators with eigenvalues
$\{ \alpha, \frac{\alpha}{4}, \frac{\alpha}{4}\}$ but the eigenspace corresponding to the multiple eigenvalue
$\frac{\alpha}{4}$ inherits a definite (positive or negative) metric in the complex case in opposition
to the paracomplex case, where the induced metric has Lorentzian signature. The latter is necessarily
the case for any Type II Jordan-Osserman metric.
\end{remark}

Authors' address

\medskip

\noindent
Department of Geometry and Topology, Faculty of Mathematics,\\
University of Santiago de Compostela, \\15782 Santiago de Compostela, Spain.

\medskip

\noindent {\it E-mails:} xtjosec@usc.es $\,\,$ xtedugr@usc.es $\,\,$ ravazlor@usc.es

\end{document}